\documentclass[a4paper,12pt]{article}

\usepackage {amssymb,latexsym}
\usepackage {amsmath}
\usepackage[french]{babel}

\newtheorem{defi}{D\'{e}finition}[section]
\newtheorem{propo}{Proposition}[section]

\newtheorem{rmk}{Remarque}

\newtheorem{exemples}{Exemples}

\begin{document}

\title{Vari\'{e}t\'{e}s de Poisson polaris\'{e}es }
\author{Azzouz AWANE\footnote{Ce travail a \'{e}t\'{e} \'{e}labor\'{e} avec l'aide de la coop\'{e}ration 
franco-marocaine Action in\'{e}tgr\'{e}e A.I. MA/02/32.} \\
Universit\'{e} Hassan II - Mohammedia. Facult\'{e} des Sciences Ben M'sik. \\
B.P. 7955. Casablanca Maroc\\
E-mail : a.awane@univh2m.ac.ma}
\date{}
\maketitle

\thispagestyle{empty} {\small \ }{\bf Abstract}. {\it We introduce and study
the basic notion of polarized Poisson manifolds generalizing the classical
case of Poisson manifolds and extend this last notion for the ${k-}$%
symplectic stuctures. And also, we show that for any polarized
Hamiltonian map, the associated Nambu's dynamical system and polarized
Hamiltonian system are connected by relations characterizing the mechanical
aspect of the $k-$symplectic geometry.} 

{\bf Keywords} : {\it Lagrangian spaces. Symplectic structure. Hamiltonian
systems. Poisson bracket.}

M.S.C. 20F05, 20F26, 51A10, 70Hxx.

\section{Introduction}

Une vari\'{e}t\'{e} polaris\'{e}e est d\'{e}finie par la donn\'{e}e sur une
vari\'{e}t\'{e} diff\'{e}rentiable $M$ de dimension paire $2n,$ d'un couple $%
\left( \theta ,E\right) $ dans lequel $\theta $ est une forme
diff\'{e}rentielle ferm\'{e}e de degr\'{e} $2$ de classe maximum et $E$ est
un sous fibr\'{e} int\'{e}grable de $TM$ de codimension $n$ annulant la $2-$%
forme $\theta ;$ en d'autres termes, $\theta $ est une structure
symplectique sur $M$ et le feuilletage ${\frak F}$ d\'{e}fini par le
sous-fibr\'{e} $E$ est lagrangien par rapport \`{a} $\theta $.

La notion de vari\'{e}t\'{e} polaris\'{e}e joue un r\^{o}le important en
th\'{e}orie de la quantification g\'{e}om\'{e}trique de Kostant-Souriau. Des
propri\'{e}t\'{e}s int\'{e}ressantes ont \'{e}t\'{e} mises en \'{e}vidence
par A.Weinstein, P.Dazord, J.M. Morvan, P. Molino, P.Libermann etc...

Le th\'{e}or\`{e}me de Darboux montre que tout point de $M$ poss\`{e}de un
voisinage ouvert muni d'un syst\`{e}me de coordonn\'{e}es locales $\left(
x^1,\ldots ,x^n,y^1,\ldots ,y^n\right) ,$ telles que 
\[
\theta =\sum_{i=1}^ndx^i\wedge dy^i 
\]
et le sous fibr\'{e} $E$ soit d\'{e}fini par les \'{e}quations $dy^1=\ldots
=dy^n=0.$

Localement, les applications hamiltoniennes polaris\'{e}es de $\left( \theta
,E\right) \ $s'\'{e}crivent : 
\[
H=\sum_{i=1}^na_i(y^1,\ldots ,y^n)x^i+b(y^1,\ldots ,y^n)
\]
et elles forment un sous-module ${\frak H}\left( M,{\frak F}\right) $ du $%
{\frak B}\left( M,{\frak F}\right) -$module ${\cal C}^\infty \left( M\right) 
$ des fonctions diff\'{e}rentiables sur $M,$ ${\frak B}\left( M,{\frak F}%
\right) $ \'{e}tant l'anneau des fonctions basiques pour le feuilletage $%
{\frak F.}$ Le tenseur de Poisson $P$ associ\'{e} \`{a} la structure
symplectique $\theta $ v\'{e}rifie de plus la relation : 
\[
P\left( dH,dK\right) \in {\frak H}\left( M,{\frak F}\right) \text{ pour tous 
}H,K\in {\frak H}\left( M,{\frak F}\right) ,
\]
ce qui nous a conduit \`{a} introduire dans ce travail la notion de
vari\'{e}t\'{e} de Poisson polaris\'{e}e sur une vari\'{e}t\'{e}
feuillet\'{e}e, permettant d'\'{e}tudier les propri\'{e}t\'{e}s de ces
nouveaux objets et de retrouver le cas usuel d'une vari\'{e}t\'{e} de
Poisson dans le cas o\`{u} ${\frak F}$ est le feuilletage trivial de
dimension $0$ dans lequel la feuille ${\frak F}_x$ passant par $x$ est
r\'{e}duite au singleton $\left\{ x\right\} $.

L'une des motivations principales qui ont conduit \`{a} introduire la
g\'{e}om\'{e}trie $k-$symplectique en tant qu'extension de la
g\'{e}om\'{e}trie de polarisation (\cite{AW-GZ}), est de proposer un support
g\'{e}om\'{e}trique des \'{e}quations de Nambu-Hamilton (\cite{NMB}), \`{a}
l'instar du formalisme hamiltonien classique, qui est une g\'{e}om\'{e}trie
de l'espace de phase (fibr\'{e} tangent $TM,$ d'une vari\'{e}t\'{e}
diff\'{e}rentiable $M,$ muni de la forme de Liouville $\lambda ).$ Rappelons
que les \'{e}quations de Hamilton 
\[
\frac{dx^i}{dt}=-\frac{\partial H}{\partial y^i},\quad \frac{dy^i}{dt}=\frac{%
\partial H}{\partial x^i}
\]
proviennent de la dualit\'{e} $X\longmapsto i(X)\theta ,$ entre les
fibr\'{e}s des rep\`{e}res et des corep\`{e}res \cite{GDB}, o\`{u} $\theta
=d\lambda ;$ et que les applications hamiltoniennes $H$ sont \`{a} valeurs
r\'{e}elles et sont reli\'{e}es aux syst\`{e}mes hamiltoniens $X_H$ par la
relation : 
\[
i\left( X_H\right) \theta =-dH.
\]

Les \'{e}quations de Nambu-Hamilton r\'{e}gissant le mouvement de la
m\'{e}canique statistique de Nambu en dimension $3$ sont donn\'{e}es par : 
\begin{equation}
\begin{array}{ll}
\frac{dx}{dt}\; & =\;\frac{D(H,G)}{D(y,z)}
\end{array}
\text{ , } 
\begin{array}{ll}
\frac{dy}{dt}\; & =\;\frac{D(H,G)}{D(z,x)}
\end{array}
\text{ , } 
\begin{array}{ll}
\frac{dz}{dt}\; & =\;\frac{D(H,G)}{D(x,y)}
\end{array}
\label{N-H}
\end{equation}
o\`{u} $\,H\,$ et $\,G\,$ sont deux fonctions r\'{e}elles d\'{e}finies sur
l'espace de phase $\,M\,$ d\'{e}crit par le syst\`{e}me de coordonn\'{e}es $%
\,(x,y,z).$

Dans cette optique, la g\'{e}om\'{e}trie $k-$symplectique propose une
structure g\'{e}om\'{e}trique dans laquelle cohabitent des $2-$formes
diff\'{e}rentielles ferm\'{e}es $\theta ^1,\ldots ,\theta ^k,$ de telle
sorte que les applications hamiltoniennes $H$ soient \`{a} valeurs dans $%
{\Bbb R}^k,$ et dont les composantes $H^p$ soient li\'{e}es au syst\`{e}mes
hamiltoniens $X_H$ par les relations : 
\[
i\left( X_H\right) \theta ^p=-dH^p, 
\]
afin de retrouver les \'{e}quations de Nambu-Hamilton tout en conservant les
traits sp\'{e}cifiques de la g\'{e}om\'{e}trie symplectique classique.

L'examen des applications hamiltoniennes d'une structure $k-$symplectique
nous ont conduit \`{a} introduire dans ce travail la notion de structure de
Poisson $k-$polaris\'{e}e sur une vari\'{e}t\'{e} feuillet\'{e}e $\left( M,%
{\frak F}\right) ,$ comme \'{e}tant un couple $\left( {\frak H}\left( M,%
{\frak F}\right) ,P\right) $ dans lequel ${\frak H}\left( M,{\frak F}\right) 
$ est un sous module de ${\cal C}^\infty \left( M,{\Bbb R}^k\right) $ sur
l'anneau des fonctions basiques ${\frak B}\left( M,{\frak F}\right) ,$ et
d'un tenseur antisym\'{e}trique 
\[
P:\bigwedge\nolimits_1\left( M,{\Bbb R}^k\right) \times
\bigwedge\nolimits_1\left( M,{\Bbb R}^k\right) \longrightarrow {\cal C}%
^\infty \left( M,{\Bbb R}^k\right) 
\]
tel que :

\begin{enumerate}
\item  pour tous $H,K\in {\frak H}\left( M,{\frak F}\right) ,$ $P\left(
dH,dK\right) \in {\frak H}\left( M,{\frak F}\right) ,$

\item  la correspondance $(H,K)\longmapsto \left\{ H,K\right\} =P\left(
dH,dK\right) ,$ de ${\frak H}\left( M,{\frak F}\right) \times {\frak H}%
\left( M,{\frak F}\right) $ \`{a} valeurs dans ${\frak H}\left( M,{\frak F}%
\right) ,$ conf\`{e}re \`{a} ${\frak H}\left( M,{\frak F}\right) $ une loi
d'alg\`{e}bre de Lie,

\item  tout \'{e}l\'{e}ment $H\in {\frak H}\left( M,{\frak F}\right) $
correspond un champ de vecteurs $X_H$ tel que $\left\langle
dK,X_H\right\rangle =\left\{ K,H\right\} .$
\end{enumerate}

Pour $k=1$, on retrouve une vari\'{e}t\'{e} de Poisson polaris\'{e}e.

Dans un syst\`{e}me de coordonn\'{e}es locales $\left( x^1,\ldots
,x^n\right) $ , le tenseur $P$ s'\'{e}crit : 
\[
P=W_{pq}^{ijr}\left( \left( \frac \partial {\partial x^i}\otimes {\omega ^p} \right) \wedge 
\left( \frac \partial {\partial x^j}\otimes \omega ^q\right)\right) \otimes e_r 
\]
o\`{u} $W_{pq}^{ijr}:U\longrightarrow {\Bbb R}$ sont des applications
diff\'{e}rentiables.

Dans la derni\`{e}re partie de ce travail, nous mettons en relief, pour
chaque \'{e}l\'{e}ment $H\in {\frak H}\left( M,{\frak F}\right) ,$ le lien
entre le syst\`{e}me hamiltonien de la structure $k-$symplectique $X_H$ et
le syst\`{e}me dynamique de Nambu $X_H^N$.

Sauf mention du contraire, les vari\'{e}t\'{e}s diff\'{e}rentiables
consid\'{e}r\'{e}es ici sont suppos\'{e}es connexes, s\'{e}par\'{e}es,
paracompactes \`{a} bases d\'{e}nombrables d'ouverts, et tous les
\'{e}l\'{e}ments introduits dans ce travail sont suppos\'{e}s de classe $%
\,C^\infty .$

\section{Vari\'{e}t\'{e}s symplectiques. Vari\'{e}t\'{e}s de Poisson}

Une vari\'{e}t\'{e} symplectique est d\'{e}finie par la donn\'{e}e d'un
couple $\left( M,\theta \right) $ dans lequel $M$ est une vari\'{e}t\'{e}
diff\'{e}rentiable de dimension paire $2n$ et $\theta $ est une $2-$forme
diff\'{e}rentielle ferm\'{e}e de classe maximum.

Le th\'{e}or\`{e}me de Darboux montre que tout point de $M$ poss\`{e}de un
voisinage ouvert $U$ muni d'un syst\`{e}me de coordonn\'{e}es locales $%
(x^1,\ldots ,x^n,y^1,\ldots ,y^n)$ tel que 
\[
\theta =\sum_{i=1}^ndx^i\wedge dy^i. 
\]

La correspondance $\zeta :X\longmapsto i\left( X\right) \theta ,$
d\'{e}finit un isomorphisme de fibr\'{e}s vectoriels au dessus de $M,$ du
fibr\'{e} des rep\`{e}res $TM$ sur le fibr\'{e} des corep\`{e}res $T^{*}M,$
permettant de d\'{e}finir un champ de bivecteurs $P$ sur $M$ par : 
\[
P\left( \alpha ,\beta \right) =-\theta \left( \zeta ^{-1}\left( \alpha
\right) ,\zeta ^{-1}\left( \beta \right) \right) 
\]
pour tous $\alpha ,\beta \in \bigwedge^1\left( M\right) .$

Un syst\`{e}me hamiltonien sur $\left( M,\theta \right) ,$ est une
transformation infinit\'{e}simale de $\theta ,$ c'est-\`{a}-dire un champ de
vecteurs $X\in {\frak X}\left( M\right) ,$ tel que $L_X\theta =0.$ Comme $%
\theta $ est ferm\'{e}e et $L_X\theta =di\left( X\right) \theta +i\left(
X\right) d\theta ,$ on d\'{e}duit qu'un syst\`{e}me hamiltonien est un champ
de vecteurs $X$ sur $M$ v\'{e}rifiant $i\left( X\right) \theta $ est
ferm\'{e}e, et donc d'apr\`{e}s le lemme de Poincar\'{e}, il existe, au
voisinage $V$ de chaque point de $M,$ une application $H:V\longrightarrow 
{\Bbb R}$, diff\'{e}rentiable telle que $i\left( X\right) \theta =-dH.$
Lorsque l'application $H$ est d\'{e}finie sur toute la vari\'{e}t\'{e} $M,$
en particulier si le premier groupe de cohomologie de de Rham est trivial,
on dira que $X$ est un syst\`{e}me hamiltonien strict et on le note par $%
X_H. $

Ainsi, \`{a} toute fonction $H$ $\in {\cal C}^\infty \left( M\right) ,$ on
peut associer, gr\^{a}ce \`{a} la dualit\'{e} $X\longmapsto i\left( X\right)
\theta ,$ un champ de vecteurs $X_H$ v\'{e}rifiant $i\left( X_H\right)
\theta =-dH,$ appel\'{e} syst\`{e}me hamiltonien associ\'{e} \`{a} $H,$ et
si, pour tout couple $\left( H,K\right) $ de fonctions diff\'{e}rentiables
sur $M,$ on pose 
\[
\left\{ H,K\right\} =-\theta \left( X_H,X_K\right) =-\theta \left( \zeta
^{-1}\left( -dH\right) ,\zeta ^{-1}\left( -dK\right) \right) =P\left(
dH,dK\right) , 
\]
l'application $\left( H,K\right) \longmapsto \left\{ H,K\right\} ,$ de $%
{\cal C}^\infty \left( M\right) \times {\cal C}^\infty \left( M\right) $
\`{a} valeurs dans ${\cal C}^\infty \left( M\right) ,$ est

\begin{enumerate}
\item  ${\Bbb R}-$bilin\'{e}aire antisym\'{e}trique,

\item  elle v\'{e}rifie l'identit\'{e} de Jacobi,

\item  elle v\'{e}rifie la formule de Leibniz : $\left\{ H,KL\right\}
=\left\{ H,K\right\} L+K\left\{ H,L\right\} .$
\end{enumerate}

La condition (3) est \'{e}quivalente \`{a} :

\begin{enumerate}
\item[(3')]  \`{a} tout \'{e}l\'{e}ment $H$ $\in {\cal C}^\infty \left(
M\right) ,$ est associ\'{e} un champ de vecteurs $X_H$ v\'{e}rifiant : 
\[
\left\langle dK,X_H\right\rangle =\left\{ K,H\right\} \text{, quel que soit }%
K\in {\cal C}^\infty \left( M\right) .
\]
\end{enumerate}

Ceci conduit \`{a} la d\'{e}finition suivante :

\begin{defi}
Une vari\'{e}t\'{e} de Poisson est un couple $\left( M,\left\{ ,\right\}
\right) $ dans lequel $M$ est une vari\'{e}t\'{e} diff\'{e}rentiable et $%
\left\{ ,\right\} $ est une application de ${\cal C}^\infty \left( M\right)
\times {\cal C}^\infty \left( M\right) $ \`{a} valeurs dans ${\cal C}^\infty
\left( M\right) ,$ v\'{e}rifiant les trois propri\'{e}t\'{e}s suivantes :

\begin{enumerate}
\item  $\left\{ ,\right\} $ est ${\Bbb R}-$bilin\'{e}aire antisym\'{e}trique,

\item  $\left\{ ,\right\} $ v\'{e}rifie l'identit\'{e} de Jacobi,

\item  \`{a} tout \'{e}l\'{e}ment $H\in {\cal C}^\infty \left( M\right) $
est associ\'{e} un champ de vecteurs $X_H$ tel que : 
\[
X_H\left( K\right) =-\left\{ H,K\right\} ,\text{ pour tout }K\in {\cal C}%
^\infty \left( M\right) .
\]
\end{enumerate}
\end{defi}

La condition 3 de cette d\'{e}finition est \'{e}quivalente \`{a} la formule
de Leibniz.

Une d\'{e}finition \'{e}quivalente est la suivante :

\begin{defi}
Une vari\'{e}t\'{e} de Poisson est un couple $\left( M,P\right) $ dans
lequel $M$ est une vari\'{e}t\'{e} diff\'{e}rentiable et $P$ est un champ de
tenseurs $2-$fois contravariant, antisym\'{e}trique, appel\'{e} bivecteur de
Poisson, tel que la correspondance 
\[
\left\{ ,\right\} :\left( H,K\right) \longmapsto \left\{ H,K\right\}
=P\left( dH,dK\right) ,
\]
de ${\cal C}^\infty \left( M\right) \times {\cal C}^\infty \left( M\right) $
dans ${\cal C}^\infty \left( M\right) $ v\'{e}rifie l'identit\'{e} de
Jacobi.
\end{defi}

\'{E}tant donn\'{e} une vari\'{e}t\'{e} de Poisson $\left( M,P\right) ,$ on
a une application lin\'{e}aire antisym\'{e}trique 
\[
\underline{P}:T^{*}M\ \longrightarrow TM, 
\]
d\'{e}finie par : 
\[
\left\langle \beta ,\underline{P}\left( \alpha \right) \right\rangle
=P\left( \alpha ,\beta \right) , 
\]
telle que pour toute fonction diff\'{e}rentiable $H:M\longrightarrow {\Bbb R}
$, le champ de vecteurs $\underline{P}(dH),$ v\'{e}rifie 
\[
\left\langle dK,\underline{P}(dH)\right\rangle =P(dH,dK)=\left\{ H,K\right\}
, 
\]
en notant $\underline{P}\left( dH\right) \,$par $-X_H,$ on obtient : 
\[
X_H(K)=-\left\langle dK,\underline{P}(dH)\right\rangle =-\left\{ H,K\right\}
=\left\{ K,H\right\} , 
\]
et donc la formule de Leibniz : 
\[
\left\{ H,KL\right\} =\left\{ H,K\right\} L+K\left\{ H,L\right\} 
\]
quelles que soient les fonctions $H,K,L\in {\cal C}^\infty \left( M\right) .$

Dans le cas o\`{u} $M$ est munie d'une structure symplectique $\theta ,$
alors localement, dans un voisinage ouvert $U$ muni des coordonn\'{e}es de
Darboux $(x^1,\ldots ,x^n,y^1,\ldots ,y^n),$ on a $\theta _{\mid
U}=dx^1\wedge dy^1+\ldots +dx^n\wedge dy^n$ et $P$ est le bivecteur sur $M$
donn\'{e} par $P\left( \alpha ,\beta \right) =-\theta \left( \zeta
^{-1}\left( \alpha \right) ,\zeta ^{-1}\left( \beta \right) \right) ,$ et
donc, 
\[
\underline{P}(dx^i)=-\frac \partial {\partial y^i}\text{ et }\underline{P}%
(dy^i)=\frac \partial {\partial x^i}, 
\]
ainsi, 
\[
P=\sum_{i=1}^n\frac \partial {\partial y^i}\wedge \frac \partial {\partial
x^i}. 
\]

\section{Syst\`{e}mes hamiltoniens polaris\'{e}s. Vari\'{e}t\'{e}s de
Poisson polaris\'{e}es}

Soit $M$ une vari\'{e}t\'{e} diff\'{e}rentiable de dimension $p+q$ munie
d'un feuilletage ${\frak F}$ de codimension $q$ et soit $E$ le sous
fibr\'{e} int\'{e}grable $p-$dimensionnel de $TM$ d\'{e}fini par les
vecteurs tangents aux feuilles de ${\frak F}\,.$ On notera par ${\Gamma }(E)$
l'ensemble des sections du $\,M-$fibr\'{e} $\,E\longrightarrow M.$

Une fonction r\'{e}elle $f$ de classe $C^\infty $ sur $M$ est dite basique
pour ${\frak F}$ si, pour tout $Y\in \Gamma (E)$, la d\'{e}riv\'{e}e $Y(f)$
de $f$ suivant $Y$ est identiquement nulle; ce qui est \'{e}quivalent \`{a}
dire que $f$ est constante sur chaque feuille de ${\frak F}$.

L'ensemble des fonctions basiques pour ${\frak F}$ sera d\'{e}sign\'{e} par $%
{\frak B}(M,{\frak F}).$ Il est clair que ${\frak B}(M,{\frak F})$ est un
sous anneau de ${\cal C}^\infty \left( M\right) $ des fonctions r\'{e}elles
diff\'{e}rentiables et d'une mani\`{e}re \'{e}vidente, ${\cal C}^\infty
\left( M\right) $ est un ${\frak B}(M,{\frak F})-$module.

Un champ de vecteurs $X\in {\frak X}\left( M\right) $ est dit feuillet\'{e}
(ou un automorphisme infinit\'{e}simal pour ${\frak F)}$ si pour tout $Y\in
\Gamma (E)$ le crochet de Lie $\left[ X,Y\right] $ appartient \`{a} $\Gamma
(E).$

Rappelons (\cite{MLN2}) que pour qu'un champ de vecteurs $X\in {\frak X}(M)$
soit feuillet\'{e}, il est n\'{e}cessaire et suffisant que si $\left(
\varphi _t\right) _{\left| t\right| <\varepsilon }$ est un groupe local
\`{a} un param\`{e}tre associ\'{e} \`{a} $X$ sur un voisinage d'un point
arbitraire de $M$, le diff\'{e}omorphisme local $\varphi _t$ laisse
invariant le sous fibr\'{e} $E,$ quel que soit $t.$

\noindent On d\'{e}signe par ${\frak L}(M,{\frak F})\,$ l'alg\`{e}bre de Lie
des champs de vecteurs feuillet\'{e}s pour ${\frak F.}$

\noindent Notons que si $f\in $ ${\frak B}(M,{\frak F})$ et $X\in {\frak L}%
(M,{\frak F}),$ la fonction $X(f)$ est basique ($X(f)\in $ ${\frak B}(M,%
{\frak F})$) et le champ de vecteurs $fX$ est un automorphisme
infinit\'{e}simal pour ${\frak F}$ ($fX\in {\frak L}(M,{\frak F})$).

Soit $M$ une vari\'{e}t\'{e} diff\'{e}rentiable de dimension paire $2n,$
munie d'une polarisation r\'{e}elle $\left( \theta ,E\right) .$ Le
th\'{e}or\`{e}me de Darboux montre que tout point de $M$ poss\`{e}de un
voisinage ouvert muni d'un syst\`{e}me de coordonn\'{e}es locales $\left(
x^1,\ldots ,x^n,y^1,\ldots ,y^n\right) ,$dites adapt\'{e}es, telles que 
\[
\theta =\sum_{i=1}^ndx^i\wedge dy^i 
\]
et le sous fibr\'{e} $E$ soit d\'{e}fini par les \'{e}quations $dy^1=\ldots
=dy^n=0.$

Un syst\`{e}me hamiltonien de la structure symplectique $\theta $ est dit
polaris\'{e}, s'il est en plus feuillet\'{e} pour le sous fibr\'{e} $E,$
autrement dit, si $X$ est une transformation infinit\'{e}simale pour la
structure symplectique $\theta $ et pour le sous fibr\'{e} $E$ \`{a} la
fois, on dira que le champ de vecteurs $X$ est un syst\`{e}me hamiltonien
polaris\'{e}$.$ Localement, d'apr\`{e}s le lemme de Poincar\'{e}, pour tout
point $x_0$ de $M,$ il existe une fonction r\'{e}elle diff\'{e}rentiable $H,$
dans un voisinage ouvert $U$ de $x_0$ telle que $i\left( X\right) \theta
=-dH.$ Par rapport \`{a} un syst\`{e}me de coordonn\'{e}es locales
adapt\'{e}es $\left( x^1,\ldots ,x^n,y^1,\ldots ,y^n\right) ,$ l'application 
$H$ et le champ de vecteurs $X$ s'\'{e}crivent : 
\[
H=\sum_{i=1}^na_i(y^1,\ldots ,y^n)x^i+b(y^1,\ldots ,y^n) 
\]
et 
\[
X=-\sum_{s=1}^n\left( \sum_{j=1}^n\frac{\partial a_j}{\partial y^s}%
(y^1,\ldots ,y^n)x^s\,\,+\,\,\frac{\partial b}{\partial y^s}(y^1,\ldots
,y^n)\right) \frac \partial {\partial x^s}+\sum_{s=1}^n\,\,a_s(y^1,\ldots
,y^n)\frac \partial {\partial y^s} 
\]
o\`{u} $a_1,\ldots ,a_n,b$ sont des fonctions basiques.

Une application diff\'{e}rentiable $H:M\longrightarrow {\Bbb R}$ est dite
hamiltonienne polaris\'{e}e s'il existe un syst\`{e}me hamiltonien
polaris\'{e} $X_H$ tel que : 
\[
i\left( X_H\right) \theta =-dH. 
\]
$X_H$ est appel\'{e} syst\`{e}me hamiltonien polaris\'{e} associ\'{e} \`{a} $%
H.$

On d\'{e}note par ${\frak H}\left( M,{\frak F}\right) $ l'ensemble des
applications hamiltoniennes polaris\'{e}es. On voit donc que, contrairement
au cas classique d'une structure symplectique, on a : 
\[
{\frak H}\left( M,{\frak F}\right) \varsubsetneq {\cal C}^\infty \left(
M\right) . 
\]
On a les propri\'{e}t\'{e}s suivantes :

\begin{enumerate}
\item  ${\frak H}\left( M,{\frak F}\right) $ est un sous ${\frak B}\left( M,%
{\frak F}\right) -$module de ${\cal C}^\infty \left( M\right) $

\item  dans un voisinage ouvert $U$ de $M,$ muni d'un syst\`{e}me de
coordonn\'{e}es locales adapt\'{e}es $\left( x^1,\ldots ,x^n,y^1,\ldots
,y^n\right) ,$ ${\frak H}\left( U,{\frak F}_U\right) $ est un sous ${\frak B}%
\left( U,{\frak F}_U\right) -$module de ${\cal C}^\infty \left( U\right) $
libre de type fini de rang $n+1,$ engendr\'{e} par $x^1,\ldots ,x^n,1.$
\end{enumerate}

Le tenseur de Poisson $P$, d\'{e}fini dans le paragraphe pr\'{e}c\'{e}dent,
v\'{e}rifie les propri\'{e}t\'{e}s suivantes :

\begin{enumerate}
\item  pour tous $H,K\in {\frak H}\left( M,{\frak F}\right) ,$ $P\left(
dH,dK\right) \in {\frak H}\left( M,{\frak F}\right) ,$

\item  la correspondance $(H,K)\longmapsto \left\{ H,K\right\} =P\left(
dH,dK\right) ,$ de ${\frak H}\left( M,{\frak F}\right) \times {\frak H}%
\left( M,{\frak F}\right) $ \`{a} valeurs dans ${\frak H}\left( M,{\frak F}%
\right) $ v\'{e}rifie la relation de Jacobi$.$

\item  \`{a} tout \'{e}l\'{e}ment $H\in {\frak H}\left( M,{\frak F}\right) $
est associ\'{e} un champ de vecteurs $X_H$ tel que : 
\[
\left\langle dK,X_H\right\rangle =\left\{ K,H\right\} 
\]
pour tout $K\in {\frak H}\left( M,{\frak F}\right) ,$ ici, on a : $X_H=-%
\underline{P}\left( dH\right) .$
\end{enumerate}

Ceci nous conduit \`{a} la notion suivante :

\begin{defi}
Soit $\left( M,{\frak F}\right) $ une vari\'{e}t\'{e} feuillet\'{e}e. On
appelle structure de Poisson polaris\'{e}e sur $M,$ un couple $\left( {\frak %
H}\left( M,{\frak F}\right) ,P\right) $ dans lequel ${\frak H}\left( M,%
{\frak F}\right) $ est un sous ${\frak B}\left( M,{\frak F}\right) -$module
de ${\cal C}^\infty \left( M\right) $ et $P:T^{*}M\times
T^{*}M\longrightarrow {\cal C}^\infty \left( M\right) $ est un bivecteur
v\'{e}rifiant les propri\'{e}t\'{e}s suivantes :

\begin{enumerate}
\item  pour tous $H,K\in {\frak H}\left( M,{\frak F}\right) ,$ $P\left(
dH,dK\right) \in {\frak H}\left( M,{\frak F}\right) ,$

\item  la correspondance $(H,K)\longmapsto \left\{ H,K\right\} =-P\left(
dH,dK\right) ,$ de ${\frak H}\left( M,{\frak F}\right) \times {\frak H}%
\left( M,{\frak F}\right) $ \`{a} valeurs dans ${\frak H}\left( M,{\frak F}%
\right) $ conf\`{e}re \`{a} ${\frak H}\left( M,{\frak F}\right) $ une loi
d'alg\`{e}bre de Lie,

\item  \`{A}\thinspace tout \'{e}l\'{e}ment $H\in {\frak H}\left( M,{\frak F}%
\right) $ correspond un champ de vecteurs $X_H$ tel que : 
\[
\left\langle dK,X_H\right\rangle =\left\{ K,H\right\} 
\]
pour tout $K\in {\frak H}\left( M,{\frak F}\right) ,$ ce champ est
d\'{e}fini par $X_H=-\underline{P}\left( dH\right) $
\end{enumerate}
\end{defi}

\begin{exemples}
\begin{enumerate}
\item  Soit $M$ une vari\'{e}t\'{e} diff\'{e}rentiable munie du feuilletage
trivial de dimension $0$ dans lequel la feuille ${\frak F}_x$ passant par $x$
est r\'{e}duite au singleton $\left\{ x\right\} $ ( ${\frak F}_x=\left\{
x\right\} )$ pour tout $x\in M.$ Dans ce cas on a :

\begin{enumerate}
\item  ${\frak B}\left( M,{\frak F}\right) ={\cal C}^\infty \left( M\right) $

\item  ${\frak H}\left( M,{\frak F}\right) $ est un ${\cal C}^\infty \left(
M\right) -$sous-module de ${\cal C}^\infty \left( M\right) ,$ donc ou bien $%
{\frak H}\left( M,{\frak F}\right) =\left( 0\right) ,$ ou bien ${\frak H}%
\left( M,{\frak F}\right) ={\cal C}^\infty \left( M\right) .$
\end{enumerate}

Deux situations se pr\'{e}sentent,

\begin{enumerate}
\item[(i)]  la structure de Poisson triviale d\'{e}finie par ${\frak H}\left( M,%
{\frak F}\right) =\left( 0\right) $ avec $P$ un bivecteur quelconque sur $M.$

\item[(ii)]  la structure classique de vari\'{e}t\'{e} de Poisson dans le cas
o\`{u} ${\frak H}\left( M,{\frak F}\right) ={\cal C}^\infty \left( M\right) $%
.
\end{enumerate}
\end{enumerate}

\item  Soit $M$ une vari\'{e}t\'{e} diff\'{e}rentiable connexe munie du
feuilletage trivial de dimension $n=\dim M.$ Dans ce cas, les fonctions
basiques sont les fonctions r\'{e}elles constantes sur $M,$ donc ${\frak B}%
\left( M,{\frak F}\right) ={\Bbb R}$, par cons\'{e}quent, ${\frak H}\left( M,%
{\frak F}\right) $ est un sous espace vectoriel du ${\Bbb R-}$espace
vectoriel ${\cal C}^\infty \left( M\right) $ et $P:T^{*}M\times
T^{*}M\longrightarrow {\cal C}^\infty \left( M\right) $ est un bivecteur tel
que ${\frak H}\left( M,{\frak F}\right) $ muni du crochet : 
\[
\left\{ H,K\right\} =P\left( dH,dK\right) 
\]
soit une sous alg\`{e}bre de Lie r\'{e}elle.
\end{exemples}
\section{Vari\'{e}t\'{e}s $\,k-$symplectiques, Vari\'{e}t\'{e}s de Poisson $%
k-$polaris\'{e}es}

\subsection{Vari\'{e}t\'{e}s $\,k-$symplectiques.}

Soit $\,M\,$ une vari\'{e}t\'{e} diff\'{e}rentiable de dimension $\,n(k+1)\,$
munie d'un feuilletage $\,{\frak F\,}$ de codimension $\,n\,$ et soient $%
\,\theta ^1,\ldots ,\theta ^k\,$ des formes diff\'{e}rentielles sur $M$
ferm\'{e}es de degr\'{e} 2.

Le sous-fibr\'{e} de $\,TM\,$ d\'{e}fini par les vecteurs tangents aux
feuilles de ${\frak F}\,$ sera d\'{e}sign\'{e} par $E$, l'ensemble des
sections du $\,M-$fibr\'{e} $E\longrightarrow M\,$ par ${\Gamma }(E).\ $

Pour tout $\,x\,$ de $\,M\,$, on d\'{e}note par $\,C_x(\theta ^1),\ldots
,C_x(\theta ^k)\,$ les sous espaces caract\'{e}ristiques des $2-$formes
diff\'{e}rentielles $\,\theta ^1,\ldots ,\theta ^k\,$ au point $\,x;$
rappelons la d\'{e}finition : 
\[
C_x(\theta ^p)=\left\{ X_x\in T_xM\text{ }\,\text{ }\mid \;i(X_x)\theta ^p=0%
\text{ et }i(X_x)d\theta ^p=0\,\right\} 
\]
o\`{u} $i(X_x)\theta ^p\,$ d\'{e}signe le produit int\'{e}rieur du vecteur $%
\,X_x\,$ par la $2-$forme $\,\theta ^p\,$.

\begin{defi}
On dit que le $\,(k+1)-$uple $\,(\theta ^1,\ldots ,\theta ^k;E)\,$ est une
structure $k-$symplectique sur $\,M\ $si pour tout $x\in M\,$, les
conditions suivantes sont satisfaites :

\begin{enumerate}
\item  $C_x(\theta ^1)\cap \cdots \cap C_x(\theta ^k)=\{0\}$,

\item  $\,\theta ^p(X,Y)\,=\,0\,$ pour tous $\,X,Y\in \Gamma (E)\,$ et $%
\,p(p=1,\ldots ,k)\,$.
\end{enumerate}
\end{defi}

Le th\'{e}or\`{e}me de Darboux montre que si $(\theta ^1,\ldots ,\theta
^k;E)\,$ est une structure $\,k-$symplectique sur $\,$la vari\'{e}t\'{e}
diff\'{e}rentiable $M\,$, alors pour tout point $x_0$ de $\,M\,$ il existe
un voisinage ouvert $\,U\,$ de $\,M\,$ contenant $\,x_0\,$ de
coordonn\'{e}es locales $\;(x^{pi},x^i)_{1\leq p\leq k,1\leq i\leq n}\;$%
dites adapt\'{e}es, tel que les formes diff\'{e}rentielles $\,\theta ^p\,$
soient repr\'{e}sent\'{e}es dans $\,U\,$ par 
\[
\theta _{\mid U}^p\,=\,\sum_{i=1}^ndx^{pi}\wedge dx^i 
\]
et le sous-fibr\'{e} $\,E\,$ soit d\'{e}fini par les \'{e}quations $%
dx^1\,=\,\ldots \,=\,dx^n\,=\,0.$

Un champ de vecteurs $\,X\,$ sur $\,M\,$ est appel\'{e} syst\`{e}me
hamiltonien $k-$polaris\'{e} (ou syst\`{e}me hamiltonien pour la structure $%
k-$symplectique) si $\,X\,$ est un automorphisme infinit\'{e}simal pour $%
{\frak F}\ $et pour les $2-$formes $\,{\theta }^p\ $\`{a} la fois; autrement
dit, s'il satisfait les conditions suivantes :

\begin{enumerate}
\item  $\,X\,$ est feuillet\'{e} pour ${\frak F}\,$,

\item  les formes de Pfaff $\,i(X){\theta }^1,\ldots ,i(X){\theta }^k\,$
sont ferm\'{e}es.
\end{enumerate}

Le champ de vecteurs $\,X\,$ sera appel\'{e} aussi automorphisme
infinit\'{e}simal pour la structure $\,k-$symplectique $\,({\theta }%
^1,\ldots ,{\theta }^k;E)$.

On d\'{e}note par ${\cal I}(M,{\cal F})$ l'espace des automorphismes
infinit\'{e}simaux pour la structure $\,k-$symplectique $\,({\theta }%
^1,\ldots ,{\theta }^k;E)$.

Le lemme de Poincar\'{e} montre que pour tout $\,x\in M\,$ il existe un
voisinage ouvert $\,U\,$ de $\,M\,$ contenant $\,x\,$ et une application
diff\'{e}rentiable $\,H\,$ de $\,U\,$ dans $\,{\Bbb R}^k\,$ dont les
composantes $H^p$ v\'{e}rifient la relation $\ $

\[
i(X){\theta }^p=-dH^p\,. 
\]
Un syst\`{e}me hamiltonien $k-$polaris\'{e} sera dit strict s'il existe une
application diff\'{e}rentiable $H:\,M\longrightarrow \,{\Bbb R}^k$ dont les
composantes $H^p$ v\'{e}rifient la relation pr\'{e}c\'{e}dente $i(X){\theta }%
^p=-dH^p\,.$ l'application $H$ est appel\'{e}e application hamiltonienne $k-$%
polaris\'{e}e et le syst\`{e}me hamiltonien $k-$polaris\'{e} $X$ sera
not\'{e} $X_H$ et sera dit associ\'{e} \`{a} $H$.

On d\'{e}signe par ${\frak H}\left( M,{\frak F}\right) $\quad le sous espace
de ${\cal C}^\infty \left( M,{\Bbb R}^k\right) $ des applications
hamiltoniennes $k-$polaris\'{e}es.

{}Soit $\,H=(H^p)_{1{\leq }p{\leq }k}\,$ une application hamiltonienne et $%
\,X_H\,$ le syst\`{e}me hamiltonien associ\'{e}. Dans un ouvert $\,U\,$ de $%
\,M\,$ muni d'un syst\`{e}me de coordonn\'{e}es locales adapt\'{e}es $%
\,(x^{pi},x^i)_{1\leq p\leq k,1\leq i\leq n}\,$, les composantes $\,H^p\,$
de $\,H\,$ et $\,X_H\,$ s'\'{e}crivent respectivement : 
\[
H^p\,=\,\sum_{j=1}^nf_j(x^1,\ldots ,x^n)x^{pj}\,+\,g^p(x^1,\ldots ,x^n) 
\]
et 
\begin{eqnarray*}
X_H\, &=&\,-\sum_{s=1}^n\sum_{p=1}^k\left( \sum_{j=1}^n\frac{\partial f_j}{%
\partial x^s}(x^1,\ldots ,x^n)x^{ps}\,\,+\,\,\frac{\partial g^p}{\partial x^s%
}(x^1,\ldots ,x^n)\right) \frac \partial {\partial x^{ps}}\,\, \\
&&+\sum_{s=1}^n\,\,f_s(x^1,\ldots ,x^n)\frac \partial {\partial x^s}
\end{eqnarray*}
o\`{u} $\,f_j\,$ et $\,g^p\,$ sont des fonctions diff\'{e}rentiables dans $%
\,U$, basiques pour le feuilletage ${\frak F}_{\mid U}$ (\cite{AW-GZ}).

Soient $\,H\,,\,K\,$ deux applications hamiltoniennes et $\,X_H\,,$ $\,X_K\,$
les syst\`{e}mes hamiltoniens associ\'{e}s. Le crochet $\,[X_H,X_K]\,$ est
un syst\`{e}me hamiltonien, et; plus pr\'{e}cis\'{e}ment, l'application
not\'{e}e $\{H,K\}$ de $\,M\,$ dans $\,{\Bbb R}^k\,$ d\'{e}finie par 
\[
\{H,K\}\,=\,-\,(\theta ^1(X_H,X_K),\ldots ,\theta ^k(X_H,X_K))
\]
satisfait \`{a} $[X_H,X_K]\,=\,X_{\{H,K\}}.${}

Dans un syst\`{e}me de coordonn\'{e}es locales adapt\'{e}es $%
\,(x^{pi},x^i)_{1\leq p\leq k,1\leq i\leq n}\,$ les composantes $\,{\{H,K\}}%
^p\,$ de $\,{\{H,K\}}\,$ s'\'{e}crivent 
\[
{\{H,K\}}^p\,=\,\sum_{s=1}^n\left( \frac{\partial H^p}{\partial x^s}\frac{%
\partial K^p}{\partial x^{ps}}\,-\,\frac{\partial H^p}{\partial x^{ps}}\frac{%
\partial K^p}{\partial x^s}\right) \,. 
\]

L'espace ${\frak H}\left( M,{\frak F}\right) $ muni du crochet $\{,\}\,$ est
une alg\`{e}bre de Lie r\'{e}elle de dimension infinie.

\subsection{Vari\'{e}t\'{e}s de Poisson $k-$polaris\'{e}es}

Soit $M$ une vari\'{e}t\'{e} diff\'{e}rentiable de dimension $n$ munie d'un
feuilletage ${\frak F.}$

Soit $\left( e_r\right) _{1\leq r\leq k}$ la base canonique de l'espace
vectoriel ${\Bbb R}^k$ dont on d\'{e}note par $\left( \omega ^r\right)
_{1\leq r\leq k}$ la base duale, et on d\'{e}note par $\bigwedge_1\left( M,%
{\Bbb R}^k\right) =\bigwedge\nolimits_1\left( M\right) \otimes {\Bbb R}^k$
l'espace des formes diff\'{e}rentielles sur $M$ \`{a} valeurs dans ${\Bbb R}%
^k,$ c'est \`{a} dire l'espace des \'{e}l\'{e}ments de la forme : 
\[
\alpha =\sum_{r=1}^k\alpha ^r\otimes e_r 
\]
avec $\alpha ^1,\ldots ,\alpha ^k$ sont des formes de Pfaff sur $M.$
Localement, dans un voisinage ouvert $U$ muni d'un syst\`{e}me de
coordonn\'{e}es locales $\left( x^1,\ldots ,x^n\right) ,$ tout
\'{e}l\'{e}ment $\alpha \in \bigwedge_1\left( M,{\Bbb R}^k\right) $
s'\'{e}crit : 
\[
\alpha _{\mid U}=\sum_{r=1}^k\sum_{i=1}^n\alpha _i^rdx^i\otimes e_r 
\]
o\`{u} $\alpha _i^r:U\longrightarrow {\Bbb R}$ sont des applications
diff\'{e}rentiables.

\begin{defi}
Soit $\left( M,{\frak F}\right) $ une vari\'{e}t\'{e} feuillet\'{e}e. On
appelle structure de Poisson $k-$polaris\'{e}e sur $M,$ un couple $\left( 
{\frak H}\left( M,{\frak F}\right) ,P\right) $ dans lequel ${\frak H}\left(
M,{\frak F}\right) $ est un sous-${\frak B}\left( M,{\frak F}\right) -$%
module de ${\cal C}^\infty \left( M,{\Bbb R}^k\right) $ et 
\[
P:\bigwedge\nolimits_1\left( M,{\Bbb R}^k\right) \times
\bigwedge\nolimits_1\left( M,{\Bbb R}^k\right) \longrightarrow {\cal C}%
^\infty \left( M,{\Bbb R}^k\right) 
\]
est une application ${\cal C}^\infty \left( M\right) -$bilin\'{e}aire
antisym\'{e}trique telle que :

\begin{enumerate}
\item  pour tous $H,K\in {\frak H}\left( M,{\frak F}\right) ,$ $P\left(
dH,dK\right) \in {\frak H}\left( M,{\frak F}\right) ,$

\item  la correspondance $(H,K)\longmapsto \left\{ H,K\right\} =P\left(
dH,dK\right) ,$ de ${\frak H}\left( M,{\frak F}\right) \times {\frak H}%
\left( M,{\frak F}\right) $ \`{a} valeurs dans ${\frak H}\left( M,{\frak F}%
\right) ,$ conf\`{e}re \`{a} ${\frak H}\left( M,{\frak F}\right) $ une loi
d'alg\`{e}bre de Lie,

\item  \`{A} tout \'{e}l\'{e}ment $H\in {\frak H}\left( M,{\frak F}\right) $
correspond un champ de vecteurs $X_H$ tel que : 
\[
\left\langle dK,X_H\right\rangle =\left\{ K,H\right\} ,
\]
pour tout $K\in {\frak H}\left( M,{\frak F}\right) .$
\end{enumerate}

$P$ sera appel\'{e} tenseur de Poisson $k-$polaris\'{e}.
\end{defi}

Par rapport \`{a} un syst\`{e}me de coordonn\'{e}es locales $\left(
x^1,\ldots ,x^n\right) $ d\'{e}fini sur un voisinage ouvert $U,$ le tenseur
de Poisson $k-$polaris\'{e} $P$ s'\'{e}crit 
\[
P=W_{pq}^{ijr}\left( \left( \frac \partial {\partial x^i}\otimes \omega ^p \right) \wedge
\left( \frac \partial {\partial x^j}\otimes \omega ^q \right)\right) \otimes e_r 
\]
o\`{u} $W_{pq}^{ijr}:U\longrightarrow {\Bbb R}$ sont des applications
diff\'{e}rentiables.

Consid\'{e}rons le cas o\`{u} $M$ est l'espace r\'{e}el ${\Bbb R}^{n\left(
k+1\right) }$ muni de la structure $k-$symplectique canonique $\left( \theta
^1,\ldots ,\theta ^k,E\right) $ d\'{e}finie par : 
\[
\theta _{\mid U}^p\,=\,\sum_{i=1}^ndx^{pi}\wedge dx^i 
\]
et le sous-fibr\'{e} $\,E\,$ soit d\'{e}fini par les \'{e}quations $%
dx^1\,=\,\ldots \,=\,dx^n\,=\,0.$

Pour tous $p=1,\ldots ,k$ et $j=1,\ldots ,n$, le syst\`{e}me hamiltonien
associ\'{e} \`{a} l'application hamiltonienne $H_{pj}=-\left( x^j\delta
^{1p},\ldots ,x^j\delta ^{qp},\ldots ,x^j\delta ^{kp}\right) $ est donn\'{e}
par : 
\[
X_{H_{pj}}=\frac \partial {\partial x^{pj}}, 
\]
et pour tout $j=1,\ldots ,n$, le syst\`{e}me hamiltonien associ\'{e} \`{a}
l'application hamiltonienne $H_j=\left( x^{1j},\ldots ,x^{kj}\right) $ est
donn\'{e} par : 
\[
X_{H_j}=\frac \partial {\partial x^j}. 
\]
Le tenseur $P$ s'\'{e}crit : 
\[
P=\sum_{p=1}^k\sum_{i=1}^n\left( \left(\frac \partial {\partial x^i}\otimes \omega
^p\right)\wedge \left( \frac \partial {\partial x^{pi}}\otimes \omega ^p \right) \right) \otimes
e_p. 
\]
On voit bien que l'on a : 
\[
P\left( dH,dK\right) =P\left( dH^q\otimes e_q,dK^r\otimes e_r\right)
=\sum_{p=1}^k\sum_{i=1}^n\left( \frac{\partial H^p}{\partial x^i}\frac{%
\partial K^p}{\partial x^{pi}}-\frac{\partial H^p}{\partial x^{pi}}\frac{%
\partial K^p}{\partial x^i}\right) e_p=\left\{ H,K\right\} 
\]
pour tous $H,K\in {\frak H}\left( M,{\frak F}\right) .$

Dans le cas o\`{u} $k=1,$ on retrouve le cas d'une vari\'{e}t\'{e}
polaris\'{e}e et le cas classique d'une vari\'{e}t\'{e} de Poisson lorsque
la vari\'{e}t\'{e} $M$ est munie du feuilletage trivial de dimension $0$
donn\'{e} par ${\frak F}_x=\left\{ x\right\} .$

Soit $\left( M,{\frak F}\right) $ une vari\'{e}t\'{e} feuillet\'{e}e de
dimension $n$ munie d'une structure de Poisson $k-$polaris\'{e}e $\left( 
{\frak H}\left( M,{\frak F}\right) ,P\right) $.

\`{A} tout \'{e}l\'{e}ment $\alpha \in \bigwedge_1\left( M,{\Bbb R}^k\right) 
$ est associ\'{e} une application ${\cal C}^\infty \left( M\right) -$%
lin\'{e}aire 
\[
P\left( \alpha ,.\right) :\bigwedge\nolimits_1\left( M,{\Bbb R}^k\right)
\longrightarrow {\cal C}^\infty \left( M,{\Bbb R}^k\right) 
\]
telle que $P\left( \alpha ,.\right) \left( \beta \right) =P\left( \alpha
,\beta \right) ,$ qui s'identifie au champ de vecteur $\underline{P}$ $%
\left( \alpha \right) $ pour $k=1.$

Pour $k$ quelconque, on a une application canonique 
\[
\Xi :{\frak X}\left( M\right) \longrightarrow {\cal L}_{{\cal C}^\infty
\left( M\right) }\left( \bigwedge\nolimits_1\left( M,{\Bbb R}^k\right) ,%
{\cal C}^\infty \left( M,{\Bbb R}^k\right) \right) 
\]
d\'{e}finie par : 
\[
\Xi \left( X\right) \left( \beta \right) =\left\langle \beta ,X\right\rangle
=\sum_{p=1}^k\beta ^p\left( X\right) e_p=\sum_{p=1}^k\left( \beta ^p\otimes
e_p\right) \left( X\right) ; 
\]
pour tous $X\in {\frak X}\left( M\right) $ et $\beta =\sum_{p=1}^k\beta
^p\otimes e_p\in \bigwedge\nolimits_1\left( M,{\Bbb R}^k\right) ;$ par
rapport \`{a} un syst\`{e}me de coordonn\'{e}es locales $\left( x^1,\ldots
,x^n\right) $ on a : 
\[
\Xi \left( \frac \partial {\partial x^l}\right) \left( \beta \right)
=\sum_{p-1}^k\frac{\partial \beta ^p}{\partial x^l}e_p. 
\]
L'application $\Xi $ est injective, et est un isomorphisme si et seulement
si $k=1.$

\begin{propo}
Soit $M$ une vari\'{e}t\'{e} diff\'{e}rentiable de dimension $n\left(
k+1\right) $ munie d'une structure $k-$symplectique $\left( \theta ^1,\ldots
,\theta ^k,E\right) $ dont on note par $\left( {\frak H}\left( M,{\frak F} \right) ,P\right)$ 
la structure de Poisson $k-$polaris\'{e}e associ\'{e}e.
Pour toute application hamiltonienne polaris\'{e}e 
$H \in {\frak H} \left( M, {\frak F} \right)$, on a $P\left( dH,.\right) \in  Im \Xi$, et plus
pr\'{e}cis\'{e}ment, le champ de vecteurs $X_H$ satisfait $:$%
\[
\Xi (X_H)=-P\left( dH,.\right) .
\]
\end{propo}

En particulier, $\left\langle dK,\Xi (X_H)\right\rangle =-P\left(
dH,dK\right) \ =\left\{ K,H\right\} $ pour tout $K\in {\frak H}\left( M,%
{\frak F}\right) .$

\begin{rmk}
Dans la d\'{e}finition d'une vari\'{e}t\'{e} $k-$polaris\'{e}e la condition
3 est \'{e}quivalente \`{a} :

\begin{enumerate}
\item[(3')]  \`{A} tout \'{e}l\'{e}ment $H\in {\frak H}\left( M,{\frak F}%
\right) $ correspond un champ de vecteurs $X_H$ tel que : 
\[
\Xi (X_H)=-P\left( dH,.\right) .
\]
\end{enumerate}
\end{rmk}

\section{Approches avec la m\'{e}canique de Nambu}

Dans son article sur la dynamique hamiltonienne g\'{e}n\'{e}ralis\'{e}e,
Y. Nambu a propos\'{e} une m\'{e}canique (\cite{NMB}), qui n'a eu jusqu'\`{a}
pr\'{e}sent que des formulations g\'{e}om\'{e}triques partielles, en dehors
de la g\'{e}om\'{e}trie $k-$symplectique qui joue un r\^{o}le central dans
le pr\'{e}sent travail, nous citons les vari\'{e}t\'{e}s de Nambu,
introduites par Sagar A. Pandit et Anil D. Gangal (\cite{PND-GGL}) qui sont
des vari\'{e}t\'{e}s diff\'{e}rentiables de dimension $3n$ munies des formes
diff\'{e}rentielles ferm\'{e}es de degr\'{e} $3$ compl\`{e}tement
antisym\'{e}trique, et strictement non d\'{e}g\'{e}n\'{e}r\'{e}e.
Localement, ces formes s'\'{e}crivent sous la forme : 
\[
\sum_{i=1}^ndx^i\wedge dy^i\wedge dz^i
\]
(Th\'{e}or\`{e}me de Nambu-Darboux), ce qui permet de retrouver les
\'{e}quations de Nambu-Hamilton (\ref{N-H}) dans le cas de $M={\Bbb R}^3$ ($%
n=1).$

Dans le cadre d'une structure $k-$symplectique, l'ensemble des applications
hamiltoniennes $k-$polaris\'{e}es ${\frak H}\left( M,{\frak F}\right) $ est
contenu strictement dans ${\cal C}^\infty \left( M,{\Bbb R}^k\right) $, mais
il se trouve que pour chaque application hamiltonienne polaris\'{e}e $H\in 
{\frak H}\left( M,{\frak F}\right) ,$ le syst\`{e}me hamiltonien de la
structure $k-$symplectique $X_H$ et le syst\`{e}me dynamique de Nambu $X_H^N$
sont li\'{e}es par les relations suivantes :

\begin{enumerate}
\item  $X_H^n\,=\left( -1\right) ^k\left( \,f(z)\right) ^{k-1}X_H$ o\`{u} $%
f\in {\frak B}\left( M,{\frak F}\right) ,$ pour $M={\Bbb R}^{k+1},$

\item  $X_H^N=\sum_{i=1}^{n} f_{i} \left( z^1,\ldots ,z^n \right) X_H^i$
o\`{u} $X_H^i=X_H \left( x^i\right) \frac {\partial} {\partial
x^i}+X_H(y^i)\frac {\partial} {\partial y^i}+X_H\left( z^i\right) \frac {\partial}{\partial z^i}$,
pour tout $i=1,\ldots ,n,$ pour  $M={\Bbb R}^{3n}.$
\end{enumerate}

Dans la premi\`{e}re \'{e}tape, on se place dans le cas o\`{u} $\,M\,$ est
l'espace r\'{e}el ${\Bbb R}^3$ muni de la structure $2-$symplectique
canonique $\,(\theta ^1,\theta ^2;E)\,$ d\'{e}finie par : 
\[
\theta ^1\;=\;dx\wedge dz,\quad \theta ^2\;=\;dy\wedge dz\text{ et }%
E\,=\,ker dz. 
\]

Les applications hamiltoniennes de la structure $\,2-$symplectique, c'est
\`{a} dire les \'{e}l\'{e}ments de ${\frak H}\left( M,{\frak F}\right) ,$
sont les applications $H\,:\,M\,\longrightarrow \,{\Bbb R}^2$ dont les
composantes sont donn\'{e}es par $H^1\;=\;f(z)x+g^1(z),$ $%
H^2\;=\;f(z)y+g^2(z),$ o\`{u} $\,f\,$, $\,g^1\,$ et $\,g^2\,$ sont des
fonctions diff\'{e}rentiables basiques d\'{e}finies sur l'espace $\,M\,$.
Les trajectoires du syst\`{e}me hamiltonien $\,X_H\,$ de la structure $2$%
-symplectique sont donn\'{e}es par les \'{e}quations suivantes : 
\[
\frac{dx}{dt}\;=\;-\frac{\partial H^1}{\partial z}\text{ , }\frac{dy}{dt}%
\;=\;-\frac{\partial H^2}{\partial z}\text{ , }\frac{dz}{dt}\;=\;\frac{%
\partial H^1}{\partial x}\;=\,\frac{\partial H^2}{\partial y} 
\]

On d\'{e}duit, que pour tout $H=(H^1,H^2)\in {\frak H}\left( M,{\frak F}%
\right) $, le syst\`{e}me hamiltonien $\,X_H\,$ et le syst\`{e}me dynamique
de Nambu $\,X_H^N\,$ sont li\'{e}s par la relation 
\[
X_H^N\,=\,f(z)X_H.
\]
et la fonction $(f(z))^{-1}H=(x+h^1(z),y+h^1(z))\,\,$ est une solution des
\'{e}quations du mouvement de la m\'{e}canique statistique de Nambu sur le
domaine de l'espace o\`{u} $f(z)$ ne s'annule pas, ici $%
h^1(z)=(f(z))^{-1}g^1(z)\,$ \quad et\qquad $\,h^2(z)=(f(z))^{-1}g^2(z)\,.$

Dans le cas o\`{u} $M$ est l'espace r\'{e}el $\ {\Bbb R}^{k+1},$ on
consid\`{e}re la structure $k-$symplectique canonique $\,(\theta ^1,\ldots
,\theta ^k;E)\,$ d\'{e}finie par : 
\[
\theta ^1\;=\;dx^1\wedge dz,\ldots ,\theta ^k\;=\;dx^k\wedge dz
\]
et $E$ est d\'{e}fini par l'\'{e}quation $dz=0,$ $\left( x^1,\ldots
,x^k,z\right) $ \'{e}tant le syst\`{e}me de coordonn\'{e}es cart\'{e}siennes
de ${\Bbb R}^{k+1}.$\quad 

Les applications hamiltoniennes de cette structure $k-$symplectique, c'est
\`{a} dire les \'{e}l\'{e}ments de ${\frak H}\left( M,{\frak F}\right) ,$
sont les applications $H\,:\,M\,\longrightarrow \,{\Bbb R}^k$ dont les
composantes sont donn\'{e}es par $H^1\;=\;f(z)x^1+g^1(z),\ldots
,H^k\;=\;f(z)x^k+g^k(z)$, o\`{u} $\,f\,$, $\,g^1,\ldots ,\,g^k\,$ sont des
fonctions basiques pour le feuilletage d\'{e}fini par $z=$cons$\tan $te,
diff\'{e}rentiables sur l'espace $\,M.$

Les trajectoires du syst\`{e}me dynamique de Nambu $X_H^N\,$ associ\'{e}
\`{a} $H$ sont donn\'{e}es par les \'{e}quations suivantes :

\[
\frac{dx^j}{dt}=\sum_{i_1,i_2,\ldots ,i_{k\ }=1}^{k+1}\varepsilon
_{ji_1i_2\ldots i_k}\frac{\partial H^1}{\partial x^{i_1}}\frac{\partial H^2}{%
\partial x^{i_2}}\ldots \frac{\partial H^k}{\partial x^{i_k}},
\]
o\`{u} $\varepsilon _{i_1i_2\ldots i_{k+1}}$ est le tenseur de Levi-Civita.
On a donc, 
\[
\frac{dz}{dt}=\varepsilon _{\left( k+1\right) 123\ldots k}\frac{\partial H^1%
}{\partial x^1}\frac{\partial H^2}{\partial x^2}\ldots \frac{\partial H^k}{%
\partial x^k}=\left( -1\right) ^k\left( f(z)\right) ^{k\ }.
\]
\[
\frac{dx^1}{dt}=\varepsilon _{1\left( k+1\right) 2\ldots k}\frac{\partial H^1%
}{\partial x}\frac{\partial H^2}{\partial x^2}\ldots \frac{\partial H^k}{%
\partial x^k}=\left( -1\right) ^{k-1}\left( \frac{\partial f(z)}{\partial z}%
x^1+\frac{\partial g^1(z)}{\partial z}\right) \left( f(z)\right) ^{k-1},
\]
\[
\frac{dx^2}{dt}=\varepsilon _{21\left( k+1\right) 2\ldots k}\frac{\partial
H^1}{\partial x^1}\frac{\partial H^2}{\partial x}\frac{\partial H^3}{%
\partial x^3}\ldots \frac{\partial H^k}{\partial x^k}=\left( -1\right)
^{k-1}\left( \frac{\partial f(z)}{\partial z}x^2+\frac{\partial g^2(z)}{%
\partial z}\right) \left( f(z)\right) ^{k-1},
\]
\[
\ldots 
\]
\[
\frac{dx^k}{dt}=\varepsilon _{k123\ldots (k+1)}\frac{\partial H^1}{\partial
x^1}\frac{\partial H^2}{\partial x^2}\ldots \frac{\partial H^{k-1}}{\partial
x^{k-1}}\frac{\partial H^k}{\partial x}=\left( -1\right) ^{k-1}\left( \frac{%
\partial f(z)}{\partial z}x^k+\frac{\partial g^k(z)}{\partial z}\right)
\left( f(z)\right) ^{k-1},
\]

Les trajectoires du syst\`{e}me hamiltonien $X_H,$ sont donn\'{e}es par :

\[
\frac{dx^1}{dt}==-\left( \frac{\partial f(z)}{\partial z}x^1+\frac{\partial
g^1(z)}{\partial z}\right) ,\ldots ,\frac{dx^k}{dt}=-\left( \frac{\partial
f(z)}{\partial z}x^k+\frac{\partial g^k(z)}{\partial z}\right) \text{ et }%
\frac{dz}{dt}=f\left( z\right) 
\]
On d\'{e}duit donc que pour tout $H=(H^1,\ldots ,H^k)$ $\in {\frak H}\left(
M,{\frak F}\right) ,$ le syst\`{e}me hamiltonien $\,X_H\,$ et le syst\`{e}me
dynamique de Nambu $\,X_H^N\,$ sont li\'{e}s par la relation 
\[
X_H^N=\left( -1\right) ^k\left( \,f(z)\right) ^{k-1}X_H. 
\]
avec $H^p=f(z)x^p+g^p(z),$ $p=1,\ldots ,k.$ Et la fonction 
\[
\left( -1\right) ^k(f(z))^{-\left( k-1\right) }H=(x^1+h^1(z),\ldots
,x^k+h^k(z)) 
\]
$\,\,$ est une solution des \'{e}quations du mouvement de la m\'{e}canique
statistique de Nambu sur le domaine de l'espace o\`{u} $f(z)$ ne s'annule
pas, ici 
\[
h^p(z)=\left( -1\right) ^k(f(z))^{-\left( k-1\right) }g^p(z)\text{ avec }%
p=1,\ldots ,k. 
\]

Enfin, pour $M={\Bbb R}^{3n},$ on munit cet espace de la structure $2-$%
symplectique canonique d\'{e}finie par : 
\[
\theta ^1=\sum_{i=1}^n dx^i \wedge dz^i,\quad \theta
^2=\sum_{i=1}^n dy^i\wedge dz^i 
\]
et $E$ est d\'{e}fini par les \'{e}quations $dz^1=\cdots =dz^n=0,$ $\left(
x^i,y^i,z^i\right) _{1\leq i\leq n}$ \'{e}tant le syst\`{e}me de
coordonn\'{e}es cart\'{e}siennes de ${\Bbb R}^{3n}.$

Dans ce cas, les applications hamiltoniennes de cette structure, c'est \`{a}
dire les \'{e}l\'{e}ments de ${\frak H}\left( M,{\frak F}\right) ,$ sont les
applications $H:M\longrightarrow {\Bbb R}^2$ dont les composantes sont
donn\'{e}es par :\ \ 
\[
H^1=\sum_{i=1}^n f_i\left( z^1,\ldots ,z^n\right) x^i+g^1\left(
z^1,\ldots ,z^n\right) \text{ et }H^2=\sum_{i=1}^n f_i\left(
z^1,\ldots ,z^n\right) y^i+g^2\left( z^1,\ldots ,z^n\right) 
\]
o\`{u} $f_1,\ldots ,f_n,$ $g^1,g^2$ sont des fonctions diff\'{e}rentiables
basiques pour le feuilletage ${\frak F}$.

Les trajectoires du syst\`{e}me dynamique de Nambu $X_H^N$ associ\'{e} \`{a} 
$H$ sont donn\'{e}es par les \'{e}quations suivantes : 
\[
\frac{dx^i}{dt}=\frac{D\left( H^1,H^2\right) }{D\left( y^i,z^i\right) }\quad
,\quad \frac{dy^i}{dt}=\frac{D\left( H^1,H^2\right) }{D\left( z^i,x^i\right) 
}\quad \text{et}\quad \frac{dz^i}{dt}=\frac{D\left( H^1,H^2\right) }{D\left(
x^i,y^i\right) }
\]
Et on v\'{e}rifie que pour toute fonction diff\'{e}rentiable $%
f:M\longrightarrow {\Bbb R}$, on a \ :

\[
\frac{df}{dt}=\sum_{i=1}^n\frac{D(f,H^1,H^2)}{D\left( x^i,y^i,z^i\right) }%
=:(H^1,H^2,f),
\]
et les trajectoires du syst\`{e}me dynamique de Nambu sont donn\'{e}es par
:\ \ \ \ \ \ \ 
\begin{eqnarray*}
\frac{dx^i}{dt} &=&\sum_{j=1}^n\frac{D\left( x^i,H^1,H^2\right) }{%
D\left( x^j,y^j,z^j\right) }=\frac{D\left( H^1,H^2\right) }{D\left(
y^i,z^i\right) }=-\left( \sum_{j=1}^n\frac{\partial f_j}{\partial z^i%
}x^j+\frac{\partial g^1}{\partial z^i}\right) f_i \\
\frac{dy^i}{dt} &=&\sum_{j=1}^n\frac{D\left( y^i,H^1,H^2\right) }{%
D\left( x^j,y^j,z^j\right) }=-\frac{D\left( H^1,H^2\right) }{D\left(
x^i,z^i\right) }=-\left( \sum_{j=1}^n\frac{\partial f_j}{\partial z^i%
}y^j+\frac{\partial g^2}{\partial z^i}\right) f_i \\
\frac{dz^i}{dt} &=&\sum_{j=1}^n\frac{D\left( z^i,H^1,H^2\right) }{%
D\left( x^j,y^j,z^j\right) }=\frac{D\left( H^1,H^2\right) }{D\left(
x^i,y^i\right) }=\left( f_i\right) ^2
\end{eqnarray*}
tandis que les trajectoires du syst\`{e}me hamiltonien $X_H$ sont
donn\'{e}es par :
\[
\begin{array}{l}
\frac{dx^i}{dt}=-\frac{\partial H^1}{\partial z^i}=-\left(
\sum_{j=1}^n\frac{\partial f_j}{\partial z^i}x^j+\frac{\partial g^1}{%
\partial z^i}\right) , \\ 
\frac{dy^i}{dt}=-\frac{\partial H^2}{\partial z^i}=-\left(
\sum_{j=1}^n\frac{\partial f_j}{\partial z^i}y^j+\frac{\partial g^2}{%
\partial z^i}\right) , \\ 
\frac{dz^i}{dt}=\frac{\partial H^1}{\partial x^i}=\frac{\partial H^2}{%
\partial y^i}=f_i.
\end{array}
\]

On d\'{e}duit que le syst\`{e}me hamiltonien polaris\'{e} $X_H$ et le
syst\`{e}me dynamique de Nambu $X_H^N$ sont li\'{e}s par la relation 
\[
X_H^N=\sum_{i=1}^n f_i \left( z^1,\ldots ,z^n \right) X_H^i 
\]

o\`{u} $X_H^i=X_H\left( x^i\right) \frac \partial {\partial
x^i}+X_H(y^i)\frac \partial {\partial y^i}+X_H\left( z^i\right) \frac
\partial {\partial z^i},$ pour tout $i=1,\ldots ,n$ .\

\end{document}